\documentclass{article}
\usepackage{xcolor}
\usepackage[linesnumbered,ruled,vlined]{algorithm2e}
\usepackage{amsthm,amsmath,amssymb}
\usepackage{pgf,tikz,pgfplots}
\usetikzlibrary{arrows}
\usepackage{mathrsfs}
\usepackage{graphicx}
\usepackage{array}
\usepackage{hyperref}
\usepackage{enumerate}
\usepackage[english]{babel}
\usepackage{longtable}
\usepackage{cite}
\usepackage{mathrsfs}
\usepackage{mathtools}

\SetCommentSty{mycommfont}
\SetKwInput{KwInput}{Input}                
\SetKwInput{KwOutput}{Output} 
\usepackage{floatrow}
\newfloatcommand{capbtabbox}{table}[][\FBwidth]
\usepackage{multirow}
\usepackage{array}
\usepackage{siunitx}
\usepackage{booktabs}

\usepackage{blindtext}
\newtheorem{theorem}{Theorem}[section]

\newtheorem{lemma}[theorem]{Lemma}

\newtheorem{dfn}[theorem]{Definition}
\newtheorem{exm}[theorem]{Example}

\title{\textbf{{An Analysis of Graceful Coloring in a Specific $r$-Regular Graphs}}}
\author{\textsc{Laavanya~ D.,  Devi Yamini~ S. }\\{\footnotesize Department of Mathematics}\\{\footnotesize Vellore Institute of Technology, Chennai, Tamil Nadu, India}\\ {\footnotesize\texttt{}}}

\date{ }
\begin{document}
\maketitle
\begin{abstract}
 A graceful $l$-coloring of a graph $G$ is a proper vertex coloring with $l$ colors which induces a proper edge coloring with at most $l-1$ colors, where the color for an edge $ab$ is the absolute difference between the colors assigned to the vertices $a$ and $b$. The graceful chromatic number $\chi_g(G)$ is the smallest $l$ for which $G$ permits graceful $l$-coloring. The problem of computing the graceful chromatic number of regular graphs is still open, though the existence of the lower bound was proved in \cite{3}. Hence, we pay attention to the computation of the graceful chromatic number of a special class of regular graphs namely complete graphs using set theoretic approach. Also, a few characterization of graphs based on their graceful chromatic number were examined.
\end{abstract}


\section{Introduction} 

Let $G=(V,E)$ be a simple connected graph with $V$ and $E$ representing the set of vertices and edges respectively. Let $|V|=n$ and $\Delta $ represent the maximum degree of the graph. Graph coloring is one of the most researched areas in graph theory. Due to the various real-time applications and conjectures in graph coloring, many mathematicians and computer scientists focused their study on this area, resulting in numerous graph coloring variants. However, all of these graph coloring variations arose from the two fundamental colorings known as vertex coloring and edge coloring. Vertex (Edge) coloring is proper if any two adjacent vertices (edges) receive distinct colors.

A graph labeling is a concept of assigning integers to the vertices or edges (or both) of a graph $G$ satisfying some conditions. It was introduced by Alexander Rosa in 1967 \cite{13}. For more details on graph labeling, refer \cite{8}. Graceful labeling, a variant of graph labeling, was initially referred as $\beta-$labeling by Alexander Rosa. A graceful labeling $g$ of a graph $G$ is a one-to-one function from the vertex set of $G$ to $\lbrace 0, 1, ..., |E(G)| \rbrace $, which induces a bijective function $g^{*}$ from $E(G)$ to $\lbrace 1, 2, ...,|E(G)| \rbrace $, defined as the absolute difference of the labels corresponding to the end vertices of the incident edge \cite{9}.

As an extension of graceful labeling, the concept of graceful coloring was introduced by Gary Chartrand \cite{3} where the vertex coloring induces the edge coloring. In graceful coloring, each color is  represented by a positive integer. For any two positive integers $i,j$ with $i < j$, let $[i,j]= \lbrace i, i+1,..., j-1, j \rbrace $. A graceful coloring of a graph $G$ is a proper vertex coloring $g $ from the vertex set of $G$ to $[1,l]$, where $l\geq 2$, which in turn induces a proper edge coloring $g ^{*}$ from the edge set of $G$ to $[1,l-1]$, such that $ g ^{*}(ab)=| g(a)- g(b)| $, for every edge $ab$ of $G$ \cite{3}. The minimum $l$ colors required for the graceful coloring of $G$ is called the graceful chromatic number of $G$, denoted by $\chi_{g}(G)$. \\

\noindent 
The following Table \ref{tab1} provides the results on graceful coloring stated in (\cite{3},\cite{6}).
\begin{longtable}[c]{|c|c|}
\hline
\textbf{Graph} \boldmath{$G$} & \begin{tabular}[c]{@{}c@{}} \boldmath {$\chi_{g}(G)$}\end{tabular}  \\ \hline

\multicolumn{1}{|l|}{\begin{tabular}[c]{@{}c@{}}Subgraph $S$ of $G$\end{tabular}}  & \multicolumn{1}{|l|}{\begin{tabular}[c]{@{}c@{}} at most $\chi_{g}(G)$\end{tabular} } \\ \hline

\multicolumn{1}{|l|}{\begin{tabular}[c]{@{}c@{}}Connected graph $G$\end{tabular}}  & \multicolumn{1}{|l|}{\begin{tabular}[c]{@{}c@{}}at least $ \Delta +1$\end{tabular} } \\ \hline

\multicolumn{1}{|l|}{\begin{tabular}[c]{@{}c@{}}Diameter at most  $2$\end{tabular}} & \multicolumn{1}{|l|}{\begin{tabular}[c]{@{}c@{}} at least $ n$\end{tabular} } \\ \hline

\multicolumn{1}{|l|}{\begin{tabular}[c]{@{}c@{}}$r-$ regular graph\end{tabular} } & \multicolumn{1}{|l|}{\begin{tabular}[c]{@{}c@{}}at least $r+2$\end{tabular} } \\ \hline

\multicolumn{1}{|l|}{\begin{tabular}[c]{@{}c@{}}Cycle $C_n$,$n\geq 4$\end{tabular}}  & $\begin{cases}
		4,& \text{if } n\neq 5 \\
		5,& \text{if } n=5
		\end{cases}$ \\ \hline
\multicolumn{1}{|l|}{\begin{tabular}[c]{@{}c@{}}Path $P_n$, $n\geq 5$\end{tabular} } & $4$  \\ \hline

  \multicolumn{1}{|l|}{ \begin{tabular}[c]{@{}c@{}}Wheel $W_n$, $n\geq 6$\end{tabular} } & $n$   \\ \hline
  
   \multicolumn{1}{|l|}{\begin{tabular}[c]{@{}c@{}}Complete bipartite \\graph of order $2$\end{tabular}}  & $n$  \\ \hline  
   
 \multicolumn{1}{|l|}{\begin{tabular}[c]{@{}c@{}}Complete tripartite \\ graph $K_{p,p,p}$ \end{tabular} } & \multicolumn{1}{|l|}{\begin{tabular}[c]{@{}c@{}}at most$\begin{cases}
		4p-1,& \text{ p is even } \\
		4p,& \text{ p is odd}
		\end{cases}$\end{tabular} } \\ \hline
 \multicolumn{1}{|l|}{\begin{tabular}[c]{@{}c@{}}$K_{3,3,3}$ \end{tabular} } & $12$  \\ \hline
 
  \multicolumn{1}{|l|}{\begin{tabular}[c]{@{}c@{}}$K_{4,4,4}$ \end{tabular} } & $15$  \\ \hline
  
    \multicolumn{1}{|l|}{\begin{tabular}[c]{@{}c@{}}Trees \end{tabular} }  & \multicolumn{1}{|l|}{\begin{tabular}[c]{@{}c@{}} at most $\lceil\frac{5\Delta}{3}\rceil$\end{tabular} } \\ \hline
    
  \multicolumn{1}{|l|}{ \begin{tabular}[c]{@{}c@{}}Caterpillars  \end{tabular}}  & \multicolumn{1}{|l|}{\begin{tabular}[c]{@{}c@{}}from $\Delta +1$ to $\Delta +2$\end{tabular}} \\ \hline
  
  \multicolumn{1}{|l|}{ \begin{tabular}[c]{@{}c@{}}Caterpillar with $d(u)=\Delta $\\ such that at most two \\$N(u)$ is of degree $\Delta $ \end{tabular}}  & $\Delta +2$ \\ \hline

\multicolumn{1}{|l|}{\begin{tabular}[c]{@{}c@{}}Rooted trees of height $2$\end{tabular} } & $\lceil \frac{1}{2}(3\Delta +1) \rceil$  \\ \hline
 \multicolumn{1}{|l|}{\begin{tabular}[c]{@{}c@{}}Rooted trees of height $3$\end{tabular}}  & $\lceil \frac{1}{8}(13\Delta +1) \rceil$ \\ \hline
  \multicolumn{1}{|l|}{\begin{tabular}[c]{@{}c@{}}Rooted trees of height $4$\end{tabular}}  & $\lceil \frac{1}{32}(53\Delta +1) \rceil$  \\ \hline
  \multicolumn{1}{|l|}{ \begin{tabular}[c]{@{}c@{}}Rooted trees with height \\at least $2+ \lfloor \frac{1}{3} \Delta \rfloor$\end{tabular} } & $\lceil \frac{5}{3} \Delta \rceil$ \\ \hline

 \caption{\label{tab1} Existing results on graceful coloring}
\end{longtable}	

There are results available on the graceful coloring of a few subclasses of unicyclic graphs \cite{1}, a few variations of ladder graphs \cite {10} and a subclass of a tree \cite{11}. The problem of computing the graceful chromatic number for bipartite graphs, complete graphs, split graphs, etc., remain open. Hence, in this paper, we compute the graceful chromatic number of a $n-1$ regular graph which is complete using set theory approach, more specifically using the concept of arithmetic progression.

\section{Graceful coloring of $K_n$} 
\begin{figure}[ht!]
	\centering
	\includegraphics[width=0.8\linewidth, height=0.10\textheight]{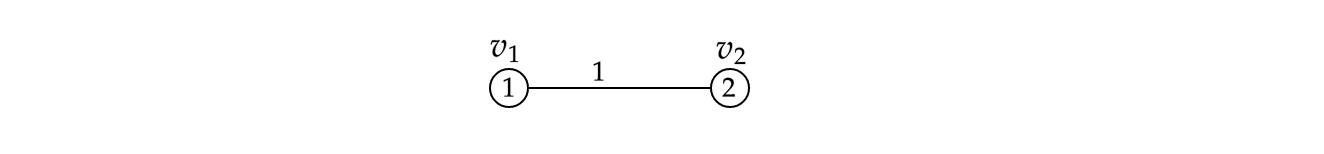}
	\caption{Graceful coloring of $K_2$}
	\label{fig1}
\end{figure}

Finding the graceful chromatic number for complete graphs $K_n$ by sketching is quite difficult, for large $n$. A complete graph $K_n$ is a graph with $n$ mutually pairwise adjacent vertices, say $v_1,v_2,\ldots,v_n$. Note that $K_2$ requires two distinct colors for the graceful coloring (refer Figure \ref{fig1}).

In order to compute $\chi_{g}(K_n)$, $n\geq 3$, we use a set theoretic approach which is discussed in the following subsection.
\subsection{$3$-AP-free sequence}
\begin{dfn}\rm \cite{2}
A series of numbers that lacks any arithmetic progression of length three is known as a $3$-AP-free sequence. Stated otherwise, a $3$-AP-free sequence lacks a triplet of distinct terms $a,b,$ and $c$ such that $|b-a|=|c-b|.$
\end{dfn}
\begin{exm}\rm
 $\lbrace 1,2,4,5,10 \rbrace$ is a $3$-AP-free sequence of length $5$. 
 \end{exm}
Given a sequence $S=\lbrace 1,2,\ldots,n\rbrace $, a subsequence $\lbrace s_1,s_2,\ldots,s_k\rbrace $ of $S$ is a $3$-AP-free sequence if it does not contain any three term arithmetic progression. Let the length of the longest $3$-AP-free sequence in $S$ be denoted by $L(n)$. Initially, Erdos and Turan in \cite{7} determined the value of $L(n)$ for $n\leq 23$ and $n=41$. Later, in \cite{14} Sharma corrected the value of $L(20)$ given in \cite{7} and also found the value of $L(n)$, for $n\leq 27$ and $n=41,42,43$. In \cite{5}, the author has found the value of $L(n)$ for $n\leq 123$ and proved the Szekeres’ conjecture for k = 5 by finding the value of $L(122)=32$. In addition, the author has given an exponential time algorithm to find the longest $3$-AP-free subsequence of $\lbrace 1,2,...,n\rbrace$ which works for all natural numbers.

\begin{lemma}\rm \label{lem1}
$K_3$ admits a graceful coloring if and only if $\lbrace g(v_1),g(v_2),g(v_3)\rbrace $ is $3$-AP-free, where $g(v_i)$ are the graceful vertex colors of $V(K_3)$.
\begin{proof}
\textbf{Necessary part:} Let $K_3$ admits a graceful coloring. Then the three mutually adjacent vertices of $K_3$ are uniquely colored such that their edge colors (absolute difference of the vertex colors) are also unique. Hence,  $\lbrace g(v_1),g(v_2),g(v_3)\rbrace $ is $3$-AP-free.\\
\textbf{Sufficient part:} Let $\lbrace g(v_1),g(v_2),g(v_3)\rbrace $ be $3$-AP-free. Then, there is no common difference between the vertex colors of $K_3$. So, $\lbrace g(v_i): i\in \lbrace 1,2,3\rbrace \rbrace$ induces a proper edge coloring in $K_3$. Hence, $K_3$ admits a graceful coloring.
\end{proof}
\end{lemma}

\begin{figure}[ht!]
	\centering
	\includegraphics[width=1.0\linewidth, height=0.14\textheight]{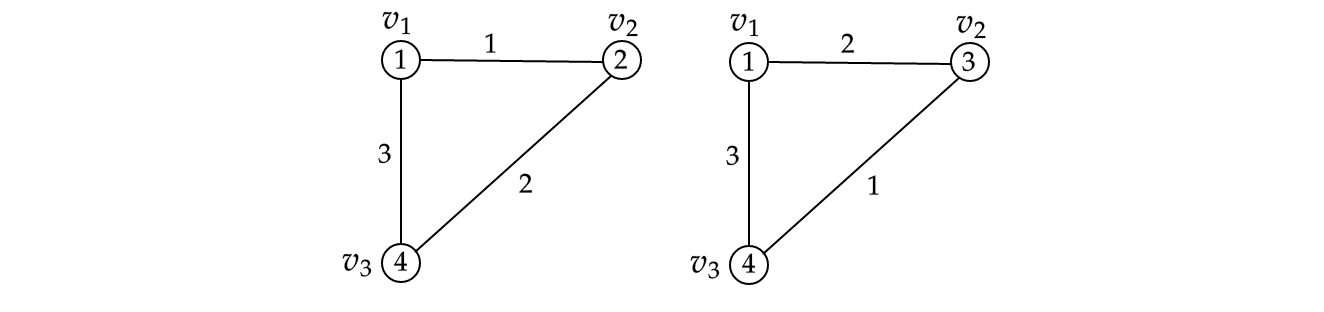}
	\caption{Graceful coloring of $K_3$}
	\label{fig2}
\end{figure}

\noindent The Figure \ref{fig2} shows two different graceful coloring of $K_3$ using four colors.

\begin{figure}[ht!]
	\centering
	\includegraphics[width=1.2\linewidth, height=0.14\textheight]{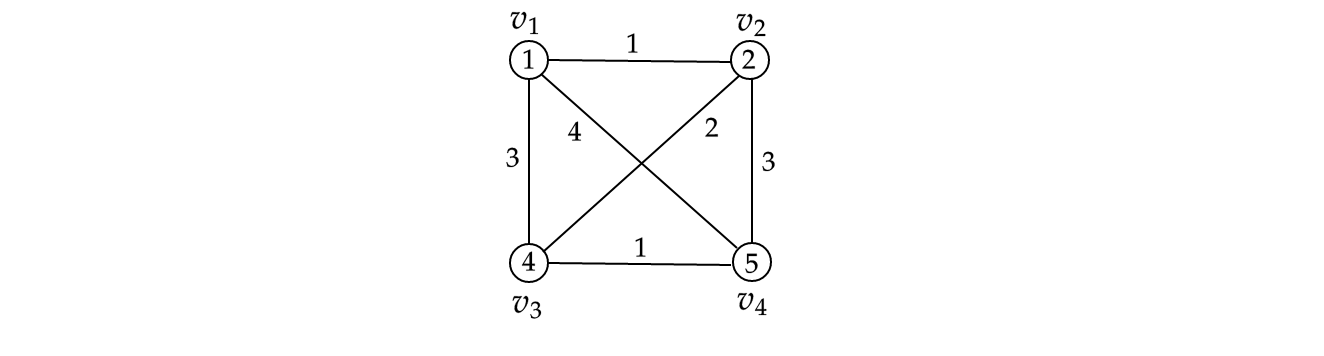}
	\caption{Graceful coloring of $K_4$}
	\label{fig3}
\end{figure}

The graceful chromatic number of $K_4$ is obtained using $g(v_1)=1$, $g(v_2)=2$, $g(v_3)=4$, and $g(v_4)=5$. Note that $\lbrace 1,2,4,5\rbrace$ is $3$-AP-free. Hence, $\chi_g(K_4)=5$. The graceful coloring of $K_4$ is represented in the Figure \ref{fig3}.

\begin{theorem}\rm
$K_n$ admits a graceful coloring if and only if $\lbrace g(v_i): v_i\in V(K_n), 1\leq i\leq n \rbrace $ is $3$-AP-free.
\begin{proof}
\textbf{Necessary part:} Suppose that $K_n$ admits a graceful coloring. Then the set of colors $\lbrace g(v_1),g(v_2),...,g(v_n)\rbrace $ of $K_n$ are proper and hence the edge colors are also proper. Since every $3$-element subset of $V(K_n)$ induces a $K_3$, which admits a graceful coloring, by the Lemma \ref{lem1}, we infer that the set of vertex colors of every $3$-element subset of $V(K_n)$ is $3$-AP-free. Hence, $\lbrace g(v_i): 1\leq i\leq n \rbrace $  is $3$-AP-free.\\
\textbf{Sufficient part:} Assume that $\lbrace g(v_1),g(v_2),\ldots,g(v_n) \rbrace $ is $3$-AP-free. Since every $3$-element subset of $V(K_n)$ induces a $K_3$, which is $3$-AP-free, by the Lemma \ref{lem1}, $K_3$ admits a graceful coloring. Hence, $K_n$ admits a graceful coloring. 
\end{proof}
\end{theorem}

 The Table \ref{tab3} describes the graceful chromatic number of $K_n$, $n\leq 32$. But $\chi_g(K_n)$ can be computed for any natural number $n$ using the algorithm given in \cite{5}.\\
\begin{longtable}{|c|c|l|}
\hline
\boldmath{$n$} & \boldmath{$\chi_{g}(K_n)$} & \begin{tabular}[c]{@{}c@{}} \boldmath{$g(v_i):1\leq i\leq n$} \end{tabular}  \\ \hline

2& 2 & {\begin{tabular}[c]{@{}c@{}}1, 2 \end{tabular} }\\ \hline

3& 4 & {\begin{tabular}[c]{@{}c@{}}1, 2, 4 \end{tabular} }\\ \hline

4& 5 & {\begin{tabular}[c]{@{}c@{}}1, 2, 4, 5 \end{tabular} }\\ \hline

5& 9 & {\begin{tabular}[c]{@{}c@{}}1, 2, 4, 8, 9\end{tabular} }\\ \hline

6& 11 & {\begin{tabular}[c]{@{}c@{}}1, 2, 4, 5, 10, 11 \end{tabular} }\\ \hline

7& 13 & {\begin{tabular}[c]{@{}c@{}}1, 2, 4, 5, 10, 11, 13 \end{tabular} }\\ \hline

8& 14 & {\begin{tabular}[c]{@{}c@{}}1, 2, 4, 5, 10, 11, 13, 14 \end{tabular} }\\ \hline

9& 20 & {\begin{tabular}[c]{@{}c@{}}1, 2, 6, 7, 9, 14, 15, 18, 20\end{tabular} }\\ \hline

10& 24 & {\begin{tabular}[c]{@{}c@{}}1, 2, 5, 7, 11, 16, 18, 19, 23, 24\end{tabular} }\\ \hline

11& 26 & {\begin{tabular}[c]{@{}c@{}}1, 2, 5, 7, 11, 16, 18, 19, 23, 24, 26\end{tabular} }\\ \hline

12& 30 & {\begin{tabular}[c]{@{}c@{}}1, 3, 4, 8, 9, 11, 20, 22, 23, 27, 28, 30\end{tabular} }\\ \hline

13& 32 & {\begin{tabular}[c]{@{}c@{}}1, 2, 4, 8, 9, 11, 19, 22, 23, 26, 28, 31, 32\end{tabular} }\\ \hline

14& 36 & {\begin{tabular}[c]{@{}c@{}}1, 2, 4, 8, 9, 13, 21, 23, 26, 27, 30, 32, 35, 36\end{tabular} }\\ \hline

15& 40 & {\begin{tabular}[c]{@{}c@{}}1, 2, 4, 5, 10, 11, 13, 14, 28, 29, 31, 32, 37, 38, 40\end{tabular} }\\ \hline

16& 41 & {\begin{tabular}[c]{@{}c@{}}1, 2, 4, 5, 10, 11, 13, 14, 28, 29, 31, 32, 37, 38, 40, 41\end{tabular} }\\ \hline

17& 51 & {\begin{tabular}[c]{@{}c@{}}1, 2, 4, 5, 10, 13, 14, 17, 31, 35, 37, 38, 40, 46, 47, 50, 51\end{tabular} }\\ \hline

18& 54 & {\begin{tabular}[c]{@{}c@{}}1, 2, 5, 6, 12, 14, 15, 17, 21, 31, 38, 39, 42, 43, 49, 51, 52, 54\end{tabular} }\\ \hline

19& 58 & {\begin{tabular}[c]{@{}l@{}}1, 2, 5, 6, 12, 14, 15, 17, 21, 31, 38, 39, 42, 43, 49, 51, 52, 54, \\58 \end{tabular} }\\ \hline

20& 63 & {\begin{tabular}[c]{@{}l@{}}1, 2, 5, 7, 11, 16, 18, 19, 24, 26, 38, 39, 42, 44, 48, 53, 55, 56, \\61, 63\end{tabular} }\\ \hline

21& 71 & {\begin{tabular}[c]{@{}l@{}}1, 2, 5, 7, 10, 17, 20, 22, 26, 31, 41, 46, 48, 49, 53, 54, 63, 64, \\68, 69, 71\end{tabular} }\\ \hline

22& 74 & {\begin{tabular}[c]{@{}l@{}}1, 2, 7, 9, 10, 14, 20, 22, 23, 25, 29, 46, 50, 52, 53, 55, 61, 65, \\66, 68, 73, 74\end{tabular} }\\ \hline

23& 82 & {\begin{tabular}[c]{@{}l@{}}1, 2, 4, 8, 9, 11, 19, 22, 23, 26, 28, 31, 49, 57, 59, 62, 63, 66, \\68, 71, 78, 81, 82\end{tabular} }\\ \hline

24& 84 & {\begin{tabular}[c]{@{}l@{}}1, 3, 4, 8, 9, 16, 18, 21, 22, 25, 30, 37, 48, 55, 60, 63, 64, 67, \\69, 76,
77, 81, 82, 84\end{tabular} }\\ \hline

25& 92 & {\begin{tabular}[c]{@{}l@{}}1, 2, 6, 8, 9, 13, 19, 21, 22, 27, 28, 39, 58, 62, 64, 67, 68, 71, \\73, 81,
83, 86, 87, 90, 92\end{tabular} }\\ \hline

26& 95 & {\begin{tabular}[c]{@{}l@{}}1, 2, 4, 5, 10, 11, 22, 23, 25, 26, 31, 32, 55, 56, 64, 65, 67, 68, \\76, 77,
82, 83, 91, 92, 94, 95\end{tabular} }\\ \hline

27& 100 & {\begin{tabular}[c]{@{}l@{}}1, 3, 6, 7, 10, 12, 20, 22, 25, 26, 29, 31, 35, 62, 66, 68, 71, 72, \\75, 77,
85, 87, 90, 91, 94, 96, 100\end{tabular} }\\ \hline

28& 104 & {\begin{tabular}[c]{@{}l@{}}1, 5, 7, 10, 11, 14, 16, 24, 26, 29, 30, 33, 35, 39, 66, 70, 72, 75, \\76, 79,
81, 89, 91, 94, 95, 98, 100, 104\end{tabular} }\\ \hline

29& 111 & {\begin{tabular}[c]{@{}l@{}}1, 2, 5, 6, 13, 15, 19, 26, 27, 30, 31, 38, 42, 44, 66, 68, 72, 77, \\80, 81,
84, 89, 93, 95, 99, 104, 107, 108, 111\end{tabular} }\\ \hline

30& 114 & {\begin{tabular}[c]{@{}l@{}}1, 2, 4, 9, 12, 13, 18, 19, 28, 30, 31, 33, 40, 45, 46, 69, 70, 75, \\82, 84,
85, 87, 96, 97, 102, 103, 106, 111, 113, 114\end{tabular} }\\ \hline

31& 121 & {\begin{tabular}[c]{@{}l@{}}1, 2, 4, 5, 10, 11, 13, 14, 28, 29, 31, 32, 37, 38, 40, 41, 82, 83, \\85, 86,
91, 92, 94, 95, 109, 110, 112, 113, 118, 119, 121\end{tabular} }\\ \hline

32& 122 & {\begin{tabular}[c]{@{}l@{}}1, 2, 4, 5, 10, 11, 13, 14, 28, 29, 31, 32, 37, 38, 40, 41, 82, 83, \\85, 86,
91, 92, 94, 95, 109, 110, 112, 113, 118, 119, 121, 122\end{tabular} }\\ \hline

\caption{\label{tab3}Graceful coloring table of $K_n$}
\end{longtable}

Equivalently, $\chi_{g}(Kn)=a(n)$, where $a(n)$ is the $n$-th entry in the integer sequence A065825 cataloged in ”The On-Line Encyclopedia of Integer Sequences” (OEIS) \cite{12}.

\section{Characterization theorems}
In this section, we give few characterization of graphs based on their graceful chromatic number. 
\begin{theorem}\rm\label{thm3}
For a simple connected graph $G$, $\chi(G)=\chi_{g}(G)$ if and only if $G\cong K_2$.
\begin{proof}
\textbf{Necessary part:} Let $\chi(G)=\chi_{g}(G)$. We need to prove that $G$ is a complete graph on two vertices. Since $\chi (G) \leq \Delta (G)+1$ (by Brooks theorem \cite{4}) and $\chi_{g}(G)\geq \Delta (G) +1$ \cite{3}, we have $\chi(G)=\chi_{g}(G)=\Delta (G)+1$. Also, $\chi (G)=\Delta (G)+1$ if and only if  the graph $G$ is an odd cycle $C_{2n+1}$, $n\geq 1$ or a complete graph. Clearly $G$ cannot be an odd cycle, as $\chi(C_{2n+1})=\Delta +1$ and $\chi_g(C_{2n+1})>\Delta +1$. Hence, $G$ is a complete graph. From the Table \ref{tab3}, we infer that $G\cong K_2$.\\
\textbf{Sufficient part:}
Suppose that $G\cong K_2$. We know that $\chi(K_2)=2=\chi_{g}(K_2)$. Hence $\chi(G)=\chi_{g}(G)$.
\end{proof}
\end{theorem}

\begin{theorem}\rm\label{thm4}
For a simple connected graph $G$, $\chi_{g}(G)=3$ if and only if $G\cong P_n$, $n\in \lbrace 3,4\rbrace$.
\begin{proof}
\textbf{Necessary part:} Assume that $\chi_{g}(G)=3$. Since $\chi_{g}(G)\geq \Delta +1$ (\cite{3}), $\Delta (G)\in \lbrace 1,2 \rbrace$. Clearly $\Delta (G)\neq 1$. If not, $G\cong K_2$ and hence $\chi_{g}(K_2)=2$, a contradiction. Hence, $\Delta (G)=2$. Observe that for $\Delta (G)=2$, $G$ is either $P_n$ or $C_n$, $n\geq 3$. From \cite{3}, $\chi_{g}(C_n)=\begin{cases}
		4,& \text{if } n\neq 5 \\
		5,& \text{if } n=5
		\end{cases} $
and $\chi_{g}(P_n)=4$ for $n\geq 5$, we conclude that $G\cong P_n$, $n\in \lbrace 3,4\rbrace$.\\
\textbf{Sufficient part:} For $G\cong P_n$, $n\in \lbrace 3,4\rbrace$, it is obvious that $\chi_{g}(G)=3$ from the following Figure \ref{fig4}.
\begin{figure}[ht!]
	\centering
	\includegraphics[width=1.0\linewidth, height=0.12\textheight]{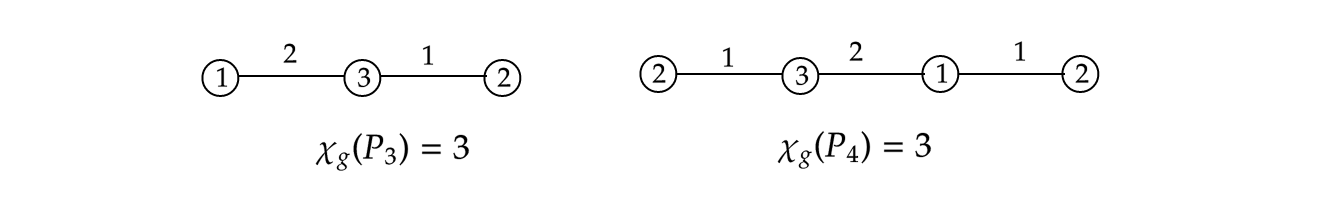}
	\caption{Graceful coloring of $P_3$ and $P_4$}
	\label{fig4}
\end{figure}
\end{proof}
\end{theorem}

\end{document}